\documentstyle[12pt,psfig]{amsart}
\def\FUN{y}

\def\figure{n}
\newtheorem{th}{Theorem}[section]	
	\newtheorem{lem}[th]{Lemma}
\newtheorem{cor}[th]{Corollary}		\newtheorem{prp}[th]{Proposition}

\newcommand{\rar}{\rightarrow} 	\newcommand{\dar}{\downarrow}
\newcommand{\lrar}{\longrightarrow}	\newcommand{\das}{\dashrightarrow}
\newcommand{\bfp}{{\Bbb{P}}}   	\newcommand{\bfc}{{\Bbb{C}}}
\newcommand{\bfq}{{\Bbb{Q}}}	
	
\newcommand{\clx}{{\cal{X}}}	
\newcommand{\co}{{\cal{O}}}	
\newcommand{\chm}{\cite{chm}}

\newcommand{\codim}{{\operatorname{codim}}}
\newcommand{\supp}{{\operatorname{supp }}}

\newcommand{\Var}{{\operatorname{Var}}}
\newcommand{\Aut}{{\operatorname{Aut}}}
\newcommand{\reld}{\operatorname{rel.dim }}

\newcommand{\red}{\mbox{\small red}}

	\newcommand{\F}{{\cal F}}
\newcommand{\G}{{\cal G}}	
\newcommand{\J}{{\cal J}}	\newcommand{\I}{{\cal I}}
	\newcommand{\OM}{\overline{\cal M}}
	\newcommand{\W}{{\cal W}}
\newcommand{\ttimes}{\tilde{\times}}

\if\FUN y
  \newcommand{\bbig}{weakly b\'{\i}g}	
  \newcommand{\subtitle}[1]{{\begin{center}{\em In Which}\\[2mm]
       \parbox{2.5in}{\small \rm{#1} }\end{center}}} 
\else 
  \newcommand{\bbig}{weakly big} 	\newcommand{\subtitle}[1]{ }
\fi

\begin{document}

\noindent %
\begin{center}
{\large  \bf
 A HIGH
FIBERED POWER OF A FAMILY OF VARIETIES OF GENERAL TYPE \\
DOMINATES A VARIETY OF
GENERAL TYPE \\[2mm]
\if\FUN y
 {\it with a few diagrams and one illustration \\ }%
\fi}

{Dan Abramovich\footnote{Partially supported by NSF grant
DMS-9503276.} \\ 
Department of Mathematics, Boston University\\
111 Cummington, Boston, MA 02215, USA \\
{\tt abrmovic@@math.bu.edu}}\\[2mm]
\today   
\end{center}

\large
\addtocounter{section}{-1}
\section{INTRODUCTION}\subtitle{We Are Introduced to Our Main
Theorem, and the Story Begins.}
\noindent %
We work over $\Bbb{C}$.

\subsection{Statement} We prove the following theorem:
\begin{th}[Fibered power theorem] \label{fibered-power}
Let $X\rar B$ be a smooth family of positive dimensional varieties of general
type, with $B$ irreducible. Then there exists an integer $n>0$, a positive
dimensional variety of 
general type $W_n$, and a dominant rational map $X^n_B \das W_n$.

Specifically, let $m_n:X^n_B \das W_n$ be the $n$-pointed birational-moduli
map. Then for sufficiently large $n$, $W_n$ is a variety of general 
type. 
\end{th}

The latter statement will be explained in section \ref{kollars}. This solves
``Conjecture H'' of \chm, \S 6.1 as well as the question at the end of
remark 1.3 in \cite{av}.

Following Viehweg's suggestions in \cite{vletter}, the fibered
power theorem is proved by way of the following theorem:
\begin{th} \label{fiberedmax}
Let $X\rar B$ be a smooth family of positive dimensional varieties of general
type, with $B$ irreducible, and $\Var(X/B) = \dim B$. Then there exists an
integer $n>0$ such that the fibered power $X^n_B$ is of general type.
\end{th}

\subsection{Main ingredients} The starting point is a theorem of Koll\'ar,
which roughly speaking says that given $f:X\rar B$ a morphism of smooth
irreducible projective varieties whose generic fiber is a variety of general
type, and $\Var(X/B) = \dim B$,  then
for large $m$ the saturation of $f_*(\omega^m_f)$ has many sections. A very
useful trick of Viehweg shows that this implies that for large $m$ the sheaf
$\omega_f^m$ itself has many sections, that is, $\omega_f$ is big. 

Following
\chm, one would like to use 
these sections pulled back to the fibered powers $f_n:X^n_B\rar B$ of $X$ over
$B$ to overcome 
any possible negativity in $\omega_B$. Unfortunately, the fibered powers may
become increasingly singular, and it is not easy to tell who wins in the race
between the positivity of $\omega_{f_n}$ and the so called adjoint conditions
imposed by 
the singularities of $X^n_B$. The fact that $\omega_f$ may have  positivity
``by accident'', as shown by the example in \chm, \S 6.1 of plane quartics,
shows that 
something more is needed - the fibration $X\rar B$ should be ``straightened
out'' before we can use sections coming from Koll\'ar's theorem. 

Semistable reduction would be sufficient for this purpose, but
unfortunately  semistable reduction for families of fiber
dimension $>2$ over a base of dimension $>1$ is yet unknown.
It is not known whether unipotent
monodromies would suffice. The case of
curves was done in  \chm\   using the moduli 
space of stable curves, and the case of surfaces was done in \cite{hassett}
using the moduli space of stable surfaces. 

Lacking such constructions in higher
dimensions, we will use a variant of semistable reduction, introduced by de
Jong \cite{dj}. This variant involves, after a suitable base change and
birational modification, a Galois cover 
$Y\rar X$, such that $Y\rar B$ is a composition of families of curves with at
most nodes. The fibers now are much better controlled, but we are left with
descending differential forms from $Y^n_B$ to $X^n_B$. This is done by studying
the behavior of the ramification ideals in the fibered powers.

\subsection{Arithmetic background and Applications} Results of this type are
motivated by Lang's conjecture. See, {\em e.g.}, \chm, \cite{hassett},
\cite{a}, \cite{av}, \cite{p}. 

Let $K$ be a number field (or any field finitely generated over $\Bbb Q$), and
let $X$ be a variety of general type over $K$. According to a well known
conjecture of S. Lang (see \cite{langbul}, conjecture 5.7), the set of
$K$-rational points $X(K)$ is not Zariski dense in $K$. In \chm, it is shown
that this conjecture of Lang implies the existence of a uniform bound $B(K,g)$
on the number of $K$ rational points on curves of genus $g$ (a stronger
implication arises if one assumes a stronger version of Lang's
conjecture). This 
result was later refined in \cite{a}, and the ultimate result of this type was
recently obtained by P. Pacelli in \cite{p}, to wit:

\begin{th}[Pacelli \cite{p}, Theorem 1.1] Assume that Lang's conjecture 
is true. Let $g\geq 2$ and $d\geq 1$ be
integers, and let $K$ be a field as above.  Then there exists an integer
$P_K(d,g)$,  such that for any
extension $L$ of $K$ of degree $d$, and any curve $C$ of genus $g$
defined over $L$, one has $$\#C(L)\leq P_K(d,g).$$
\end{th}

The main geometric ingredient in the above mentioned results is the
``Correlation Theorem'' of \chm, which is Theorem \ref{fibered-power} for
curves. In \cite{chm}, \S 6, a version of Theorem \ref{fibered-power} was
conjectured (``Conjecture H''), with the suggestion that strong uniformity
results would follow from such a theorem. This was further investigated in
\cite{av}, where it is shown that Theorem \ref{fibered-power} gives an
alternative proof for Pacelli's
theorem, as well as other strong implication results for curves and higher
dimensional varieties. As an example of a result which does not follow from
Pacelli's theorem, we have (see \cite{av}, Corollary 3.4 and theorem 1.7):

\begin{cor} Assume that the weak Lang conjecture holds.  Fix a number field $K$
and an integer $d$. Then there is a uniform bound $N$ for the number of points
of degree $d$ over $K$ on any curve $C$ of genus $g$ and gonality $>2d$ over
$K$. In fact, $N$ depends only on $g,d$ and the degree $[K:\bfq]$.  
\end{cor}

Another conjecture of Lang states that if $X$ is a complex variety of general
type, the union of rational curves on $X$ is not Zariski dense. In
\cite{av} a few implications of this conjecture were investigated, see
\cite{av} \S 3. Using \cite{p}, proposition 2.8, P. Pacelli is able to
obtain the following remarkable result (\cite{p}, corollary 5.4): Lang's
conjecture about rational curves implies that there is a bound $P_{geom}(d,g)$
such 
that for any complex curve $C$ of genus $g$, and curve $D$ of gonality $d$, the
number of nonconstant morphisms $D\rar C$ is bounded by $P_{geom}(d,g)$. This
result can again be deduced from theorem \ref{fibered-power}, but in an
unsatisfactory way - one has to reprove proposition 2.8 of \cite{p}, repeating
some of the steps, and therefore Pacelli's direct method is more appropriate. 

Another direction where our theorem falls short of obtaining a definite result
is the logarithmic case. Here again Pacelli's methods should directly imply
results regarding uniform boundedness of stably integral points on elliptic
curves (see \cite{a-int} for the definition). One suspects that in the future a
fibered power theorem will be available for log-varieties. At the moment, the
main difficulties seem to lie in obtaining an appropriate generalization of the
theorems of Koll\'ar and Viehweg.

\subsection{Acknowledgements} I would like to thank  F. Hajir and R. Gross,
whose question kept me thinking about the problem through a period when no
significant result seemed to be possible. Thanks to B. Hassett, J. de Jong,
J. Koll\'ar, P. Pacelli, R. Pandharipande, E. Viehweg  and F. Voloch, 
for helpful discussions. The realization that pluri-nodal fibration have
mild singularities was inspired by Hassett's results in \cite{hassett}, \S 4.
The understanding of the utility of such fibrations was reinforced by Pacelli's
results. 

\section{PRELIMINARIES}\subtitle{We Set Up Some Terminology
about Families of Varieties and State a Lemma about Groups.}
\subsection{Definitions}

 A variety is called an {\bf r-G variety} if it has only rational
Gorenstein (and hence canonical) singularities. For a Gorenstein variety $X$ to
be r-G, it is necessary and sufficient that for any resolution of singularities
$r:Y \rar X$ one has $r_*\omega_Y = \omega_X$ (see \cite{elkik}, II).

We say that a flat morphism of irreducible varieties $Y\rar B$ is {\bf mild},
if for any dominant $B_1 \rar B$ where $B_1$ is r-G, we have that ${Y_1}  =
Y\times_B B_1$ is r-G. Note that, by induction, if  $Y\rar B$ is mild then the
fibered powers $Y^n_B\rar B$ are mild as well.

An {\bf alteration} is a projective, surjective, generically finite morphism 
of irreducible varieties. This is slightly different from the definition in
\cite{dj}, where propernes is assumed rather than projectivity. An
alteration $B_1\rar B$ is {\bf Galois} if there 
exists a finite group  $G\subset\Aut_B B_1$ such that $B_1/G\rar B$ is
birational. 

A {\bf fibration} is a projective morphism of irreducible normal varieties
whose 
general fibers 
are irreducible and normal.

Given a fibration $X\rar B$ and an alteration $B_1\rar B$ we denote by
$X\ttimes_B B_1$ the {\bf main component } of $X\times_B B_1$. Namely, if
$\eta_{B_1}$ is the generic point of $B_1$, then $X\ttimes_B B_1 =
\overline{X\times_B \eta_{B_1}}$.

A {\bf family} is a flat fibration.

A family $Y\rar Y_1$ is called a {\bf nodal fibration} if every fiber is a
curve with at most ordinary 
nodes. A family $Y\rar B$ is called a {\bf pluri-nodal fibration} if it is a
composition of nodal fibrations $Y\rar Y_1\rar \cdots \rar B$. Note that while
nodal fibrations are generically smooth, this is not the case with pluri-nodal
fibrations. 

Given a line bundle $L$ and an ideal sheaf $\I$ on a variety $X$, we say that
$L\otimes \I$ is {\bf big} if for some $k>0$ the rational map associated to 
 $H^0(X, L^{\otimes k}\otimes \I^k)$ is birational to the image. It readily
follows that, if $L\otimes \I$ is big, and $J$ is an ideal sheaf, then for
sufficiently large $k$ we have that $L^{\otimes k}\otimes \I^kJ$ is big. The
definition differs 
somewhat from Koll\'ar's definition in \cite{kollar}. In section \ref{kollars}
we will refer to sheaves which 
are ``big''  in Koll\'ar's sense as {\em \bbig}: we say that a sheaf $\F$ is
{\bf \bbig}, if for any ample $L$ there is an integer $a$ such that
$Sym^a(\F)\otimes L^{\otimes (-1)} $ is weakly positive (see \cite{kollar},
p.367, (vii)). 

\subsection{Group theory} For a finite group $G$ let $e(G) =
\operatorname{lcm}\{\operatorname{ord}(g)|g\in G\}$. We will use the following
obvious lemma:
\begin{lem}
Let $G$ be a finite group. Then for any $n$, we have $e(G^n) = e(G)$.
\end{lem}
\section{RAMIFICATION}\subtitle{We Encounter Ramification Ideals
Measuring Differences Between Pluricanonical Sheaves in a Quotient Situation,
and Show that These Ideals Can Be Controlled in Various Cases.}
Let $V$ be a quasi projective r-G variety, $G\subset \Aut(V)$ a finite
group. Let $W = 
V/G$, and $q:V\rar W$ the quotient map. Let $r:W_1 \rar W$ be a
resolution of singularities. Note that $W$ is normal, therefore it is regular
in codimension 1.  We can pull back sections of pluricanonical sheaves on
the nonsingular locus $W_{\mbox{\small ns}}$ and extend them into the
pluricanonical sheaf of $V$. Thus, for an 
integer 
$n>0$ we have an injective morphism $\phi_n: q^*r_*\omega^n_{W_1} \rar
\omega^n_V$.  

Define the {\bf $n$-th ramification ideal} $\J_n =\J_n(G,V)=
\operatorname{Ann}\operatorname{Coker} \phi_n$.  

\begin{lem}\begin{enumerate} \item We have $\J_n\otimes \omega^n_V \cong
q^*r_*\omega^n_{W_1}$.
\item For any integer $k>0$ we have $\J_n^k \subset \J_{nk}$.
\item The ideals $\J_n$ are locally defined: if
$V'\subset V$ is a $G$-invariant open subset, then
$\J_n(G,V')=\J_n(G,V)|_{V'}$. 
\item The ideals $\J_n$ are independent of the choice of resolution $r:W_1 \rar
W$. 
\item The ideals  $\J_n$ can be also obtained if we use a partial resolution
$r:W_1 \rar W$ where $W_1$ is r-G.
\end{enumerate}
\end{lem}

{\bf Proof.} 
Since $\omega_V$ is by assumption
invertible, we have (1). 
For the same reason (2) follows: if
$\omega=\prod_{i=1}^k \omega_i$ where $\omega_i$ are local sections of
$\omega_V^n$ and if $f=\prod_{i=1}^k f_i$ where $f_i\in \J_n$, then we can
write $f_i\cdot\omega_i = \sum g_{i,j}\cdot(q^*r_* \eta_{i,j})$ and expanding
we get that $f\omega$ is a local section of  $q^*r_*\omega^{nk}_{W_1}$. It
would be nice to have an actual equality for high $n$. Part (3) follows by
definition. Parts (4) and (5) follow
by noticing that for a birational morphism $r':W_2\rar W_1$ with $W_2$
nonsingular, we have $r'_*\omega_{W_2}^n =\omega_{W_1}^n$ in both cases. \qed

Traditionally, one studies ramification by reducing to the case where both $V$
and $W$ are regular. Most of the results below follow that line, with the
exception of Proposition \ref{product}, where the author finds it liberating,
if not essential, to avoid unnecessary blowups. 

The ideals $\J_n$ give conditions for invariant differential forms
to descend to regular forms on the quotient:

\begin{prp}\label{invariants} Given an integer $n>0$ we have 
$$ (q_* (\omega^{n}_V\otimes \J_n))^G = r_* \omega^{n}_{W_1}.$$
\end{prp}
{\bf Proof:} A local section of $ (q_* (\omega^{n}_V\otimes \J_n))^G$ can be 
written 
as $\sum q_*(f_i) r_*(s_i)$, where $f_i$ are $G$ invariant, therefore $f_i =
q^*g_i$. \qed

The ideals $\J_n$ are bounded below in terms of multiplicities (here we first
use the assumption on rational singularities):

\begin{prp}[\chm\ \S 4.2, lemma 4.1]\label{multipl} Let $\Sigma_{G,V}=\Sigma
\subset V$ be the locus of fixed points: 
$$ \Sigma = \{x\in V | \exists g\in G, g(x) = x\},$$ viewed as a closed reduced
subscheme, with ideal $\I_\Sigma$. Then $\I_{\Sigma}^{n\cdot (e(G)-1)}\subset
\J_n$. 
\end{prp} 
{ {\bf Proof.} Let $V_1$ be the normalization of $W_1$ in ${\Bbb
C}(V)$. Let $W_1'\subset W_1$ be the open subset over which both $V_1$ and the
branch locus $B_{V_1/W_1}$ are nonsingular. The codimension of $W_1\setminus
W_1'$ is at least 2. Let $V_1'$ be the inverse image of $W_1'$. We have a
diagram 
$$\begin{array}{lcl} V_1' & \stackrel{s}{\rar} & V \\
                     \dar q_1 &            & \dar q \\
                     W_1' & \stackrel{r'}{\rar} & W 
 \end{array}
$$
 Let $\omega$ be a $G$-invariant $n$-canonical form on $V$, vanishing to order
$n\cdot 
(e(G)-1)$ on $\Sigma$. To show that $\omega$ descends to $W_1$ it suffices to
descend it to $W_1'$, since the codimension of the complement is at least
2. Since $V$ is r-G,  $\omega' = s^*\omega$ is a regular $n$-canonical form on
$V_1'$, vanishing to 
order $n\cdot (e(G)-1)$ on $B_{V_1/W_1}$. The 
subgroup fixing a general point of a component of $B_{V_1/W_1}$ is cyclic, and
the action is given formally by $u_1\mapsto \zeta_k u_1, u_i\mapsto u_i$ for
some root of unity 
$\zeta_k, k\leq e(G)$, where $u_i$ are local parameters, $u_1$ a uniformizer
for $B_{V_1/W_1}$. Formally at such a point, the quotient 
map is given by $w_1=u_1^k, w_i=u_i, \, i>1$. By assumption, $\omega'$ can be
written 
in terms of local parameters as $\omega'=f(u) (u^{k-1}du_1\wedge\cdots\wedge
du_m)^n=f(u) {q_1}^*(dw_1\wedge\cdots\wedge dw_m)^n$. The invariance implies
that $f(u) = {q_1}^*g(w)$ and therefore  $\omega' = {q_1}^* g(w)
(dw_1\wedge\cdots\wedge dw_m)^n$.}\qed

{\bf Remark.} It is not difficult to obtain the following refinement of this
proposition (see analogous case in \cite{kollar}, lemma 3.2): let
$B=q(\Sigma)_{\red}$, and let $I_B$ be the defining ideal. Then
$q^{-1}\I_B^{\lfloor n(1-{1\over{e(G)}})\rfloor}\subset J_n$.   

Recall that if a group $G$ acts on a variety $V$, a line bundle $L$ and an
ideal $\I$ then the ring of invariant sections $\oplus_{k\geq 0}
H^0(Y,L^{\otimes k}\otimes 
\I^k)^G$ has the same dimension as the ring of sections $\oplus_{k\geq 0}
H^0(Y,L^{\otimes k}\otimes \I^k)$. This allows us to have:

\begin{cor}[See more general statement in 
\cite{p}, Lemma 4.2]\label{prod-quotient} 
Let $X$ be a variety of general type and let 
$G=\Aut_\bfc(\bfc(X))$ be its birational automorphism group. Then for some
$n>0$ 
the quotient variety $X^n/G$, where $G$ acts diagonally, is of general type.
\end{cor} 
{
{\bf Proof.} Applying Hironaka's equvariant resolution of singularities, we may
assume that $X$ is regular and $G=\Aut X$. Let $p_i:X^n\rar X$ be
the projection onto the $i$-th factor. Choose $n$ large enough
so that $\omega^n_X\otimes \I_{\Sigma_{G,X}}$ is big. Therefore
$\omega^n_{X^n}\otimes (\sum p_i^{-1} \I_{\Sigma_{G,X}})^n$ is big. But $$(\sum
p_i^{-1}\I_{\Sigma_{G,X}})^n \subset \I_{\Sigma_{G,X^n}} \subset
J_n(G,X^n),$$ giving the result.
}\qed

Let $\Sigma \subset V$ be the locus of fixed points, and let
$\Sigma = \Sigma_1 \cup \Sigma_2 $ be a closed decomposition. Then $\J_n$ is
supported along $\Sigma$, and can be written as
$\J_n=\J_{n,\Sigma_1}\cap \J_{n,\Sigma_2}$.
Applying \ref{multipl} we obtain:
\begin{cor}\label{decompose-ramification} We have $
(\I_{\Sigma_2}^{e(G)})^n \cdot\J_{n,\Sigma_1}\subset \J_n$. 
\end{cor}

Our goal is to apply our propositions to powers of mild families.
 First, let $f:V\rar B$ be mild. Assume that $B$ is
r-G. As before, let 
$G\subset 
\Aut_B(V)$, $W = V/G$, and 
$q:V\rar W$ the quotient map. 

Let $p_i: V^m_B\rar V$ be the $i$-th projection. 
We naturally have $G^m\subset \Aut_B(V^m_B)$ acting on all components.
We denote by $q_m:V^m_B \rar W^m_B$ the associated map. Let $r:W_1\rar W$ be a
resolution of singularities. 

Define $\J_{m,n} = \prod p_i^{-1} \J_{n}$. 

\begin{lem}\label{smoothproduct} Assume that $W_1\rar B$ is mild. Then
 $\J_{m,n} \subset \J_n(G^m,V^m_B)$.
\end{lem}

{\bf Proof.} Denote $r_m:W_m = (W_1)_B^m\rar W^m_B  $ and $p_{i,W}:W_m\rar W_1$
the $i$-th projection. Since $V\rar B$ and
$W_1\rar B$ are mild, we have that 
$$\omega^n_{V^m_B/B} = 
\otimes_i \, p_i^* (\omega^n_{V/B}) \quad \mbox{ and }\quad \omega^n_{W_m/B}
= \otimes_i \, p_{i,W}^* (\omega^n_{W_1/B}).$$ Suppose a local section $w$ of
$\omega^n_{V^m_B/B}$ is 
a monomial written as $w = \prod p_i^*w_i$, 
and suppose $f\in \J_{m,n}$ is a monomial written as $f=\prod p_i^*f_i$. Then
$fw = \prod  p_i^*(f_i w_i)$ is a local section of
$q_m^*{r_m}_*\omega^n_{W_m/B}$. \qed

\begin{prp}\label{product} There exists a closed subset $F\subset B$ such that
 $$(\I_F^{e(G)})^n\cdot\J_{m,n} \subset \J_n(G^m,V^m_B).$$
\end{prp}

{\bf Proof.} Let $F\subset B$ be the discriminant locus of $W_1\rar B$, and
$U=B\setminus F$. Now apply \ref{smoothproduct} and 
\ref{decompose-ramification}. \qed

{\bf Remark.} It follows from the remark after \ref{multipl} that already
$$(\I_F^{\lfloor n(1-{1\over{e(G)}})\rfloor})\cdot\J_{m,n} \subset
\J_n(G^m,V^m_B).$$ 

We will often need to perform base
changes for fibrations. We need to find a condition on the base changed
fibration which guarantees that the original variety is of general type. This
is provided by the following proposition (which is probably well known):
 
\begin{prp}\label{base-change} Given an alteration $\rho:B_1 \rar B$ between
smooth projective
varieties,
there exists an ideal sheaf 
$\I\subset \co_{{B_1}}$ with the following property: given a fibration $f:Y\rar
B$, with ${Y_1}\rar 
Y\ttimes_B B_1$ a resolution of singularities, 
${f_1}: {Y_1} \rar {B_1}$ the induced projection,
such that 
$\omega_{{f_1}} \otimes {f_1}^{-1}\I$ is big, then $Y$ is of general
type.
\end{prp}

First a lemma:
\begin{lem} \begin{enumerate} 
\item Let $g:Y_1\rar Y$ be a generically finite morphism of smooth projective
varieties. Let $B\subset Y$ be the branch locus. Then
there exists an effective $g$-exceptional divisor $E$ on $Y_1$ and an injection
$\omega_{Y_1}(-g^*B) \rar g^*\omega_Y\otimes \co_{Y_1}(E)$.
\item If $\omega_{Y_1}(-g^*B)$ is big, then $\omega_Y$ is big as well.
\end{enumerate}
\end{lem} 

{\bf Proof.}  The pull-back morphism $g^*\omega_Y\rar \omega_{Y_1}$  gives
$g^*\omega_Y=  \omega_{Y_1}(-R-E)$ where $E$ is an effective exceptional
divisor and $R$ is the ramification divisor. Clearly $R < g^*B$.

Assume that   $\omega_{Y_1}(-g^*B)$ is big. Then $ g^*\omega_Y\otimes
\co_{Y_1}(E)$ is big. Let $ Y_1\stackrel{g_1}{\rar} Y'\stackrel{s}{\rar} Y$ be
the Stein factorization. Since $Y'$ is normal and $E$ is $g_1$-exceptional we
have that $s^*\omega_Y\otimes {g_1}_*\co_{Y_1}(E)=s^*\omega_Y$ 
therefore 
$s^*\omega_Y$ is big. Since $s$ is finite we have that $\omega_Y$ is big.

{\bf Proof of \ref{base-change}.} Choose a nonzero ideal $\I_0\subset
\co_{B_1}$ with an injection $\I_0\subset \omega_{B_1}$, and an ideal
$\I_1\subset \co_{B}$ such that $\omega_{B_1}\otimes\rho^{-1}\I_1\subset
\rho^*\omega_B$. Given a fibration $f:Y\rar B$, with $g:Y_1\rar Y$ as above, we
have that the ideal $\I_1$ vanishes on the branch locus of $g$. Set $\I =
\I_1\rho^{-1}\I_2$.  Assume that $\omega_{Y_1/B_1}\otimes g^{-1}\I$ is big,
then 
$\omega_{Y_1}\otimes (\rho\circ g)^{-1}\I_2$ is big, therefore
$\omega_{Y_1}(-g^*B)$ is big, and by the lemma we have that $\omega_Y$ is
big. \qed

\section{MAXIMAL VARIATION AND KOLL\'AR'S
THEOREMS}\label{kollars}\subtitle{We Reduce Our Theorem to the Maximal
Variation Case, and Quote a Big Theorem Of Koll\'ar Producing Sections.}

\subsection{Pointed birational moduli} The following is an
immediate generalization of Koll\'ar's generic moduli theorem (\cite{kollar},
2.4):

\begin{th}[Pointed birational moduli theorem] Let $f:X\rar B$ be a family of
varieties of general type. There exist open sets
$U\subset B$ and $V\subset f^{-1}U$, and projective varieties $Z$ and
$W_n,\quad n\geq 1$, with a diagram:
$$\begin{array}{lclclcl}
 V^n_B   & \rar & V^{n-1}_B      & \rar & \cdots & \rar & U \\
 \dar m_n &       & \dar m_{n-1} &      &        &      &\dar  m_0 \\
 W_n     & \rar   & W_{n-1}     & \rar & \cdots & \rar & Z 
\end{array}
$$
satisfying the following requirements:
\begin{enumerate}
\item The morphisms $m_n$ are dominant.
\item If $P=(P_1,\ldots,P_n), P'=(P'_1,\ldots,P'_n)\in V^n_B, f_n(P)=b,
f_n(P')=b'\in U$, then 
$m_n(P) = m_n(P')$ if and only if there exists a birational map $g:V_b\das
V_{b'}$ which is defined and invertible at $P_i$, such that $g(P_i) = P'_i$. 
\item For general $b\in U$, let $G$ be the birational automorphism group of
$X_b$, then the fiber of $W_n$ over $m_0(b)$ is birational to $X_b^n/G$, where
$G$ acts diagonally.
\item There are canonical generically finite rational maps
$W_{nk}\das (W_n)^k_Z$. 
\end{enumerate}
\end{th}

{\bf Sketch of proof:} Parts (3) and (4) follow from (2). The proof of (1) and
(2) is a 
simple modification of \cite{kollar}, 2.4, where we let $PGL$ act on the
universal family over the Hilbert scheme and its fibered powers.\qed

\subsection{Reduction of theorem \ref{fibered-power} to theorem
\ref{fiberedmax}} Recall 
by corollary \ref{prod-quotient} that for sufficiently large $n$ the 
general fiber of $W_n\rar Z$ is of general type. Also,  a simple lemma below
shows that for large $n$ the
family $W_n\rar Z$ is of maximal variation.
Assuming that theorem \ref{fiberedmax} holds true, we have that for large $k$
the variety $(W_n)^k_Z$ is of general type, therefore $W_{nk}$ is of general
type. For any $n'>nk$, applying the additivity theorem (Satz III of
\cite{viehweg1})  to
$W_{n'}\das W_{nk}$  we have that 
 the variety $W_{n'}$ is of general type. Therefore $X^{n'}_B$
dominates a variety of general type.

\begin{lem} Suppose $X\rar B$ is a one dimensional family of varieties of
general type, $\Var(X/B) =1 $, and $G\subset\Aut_BX$.  Then for sufficiently
large 
$n$, the quotient by the diagonal action $ W_n = X^n_B/G \rar B$ has
$\Var(W_n/B)=1$. 
\end{lem}

{\bf Proof.} This is immediate from the theorems of Kobayashi-Ochiai (see
\cite{d-m})  and Maehara (see \cite{moriwaki}). { Using Proposition
\ref{prod-quotient}, choose $n$ so that
the general fiber of $W_{n_0}$ over $B$ is of general type. We show that
$Var(W_{n_0+1}/B)=1$, and by induction this follows for any higher $n$. Assume
the opposite. We have the projection map $W_{n_0+1}\rar W_{n_0}$. The theorem
of Maehara implies that $Var(W_{n_0}/B)=0$: a family of varieties of genral
type dominated by a fixed variety is isotrivial. The theorem of
Kobayashi-Ochiai implies that the map  $W_{n_0+1}\rar W_{n_0}$ is isotrivial: a
family of rational maps from a fixed variety to a fixed variety of general type
is isotrivial. But the general fiber of $(W_{n_0+1})_b\rar (W_{n_0})_b$ is
isomorphic to $X_b$ (one only needs to avoid the fixed point set of the
action) - implying that $X\rar B$ is isotrivial.}\qed

\subsection{Koll\'ar's big theorem} Here we introduce the main source
for global sections.

\begin{th}[Koll\'ar's big theorem, \cite{kollar}, I, p. 363] 
 Suppose that $\pi:X\rar B$ is a fibration of positive dimensional varieties of
general type, and $\Var(X/B) = \dim B$. Assume both $X$
and $B$ are smooth. 
There is an  integer $n>0$ such that 
the sheaf $\pi_*\omega_\pi^n$ is {\em \bbig}. \qed
\end{th}

Koll\'ar's use of {\em \bbig} requires saturations, which means that the
sections obtained may have poles over exceptional divisors of the map $X\rar
B$. From this one first deduces:

\begin{cor}[\cite{viehweg}, Corollary 7.2]  Suppose that $\pi:X\rar B$ is as
above. There is a divisor $D$ on $X$ such that $\codim (\pi(\supp 
B))>1$, and  such that $\omega_\pi(D)$ is big. \qed
\end{cor}

We still have the annoying divisor $D$. Our method below will allow us to
ignore it, but actually a trick of Viehweg (\cite{viehweg}, lemma 7.3) makes it
easier. Viehweg simply applies the theorem above to $X'\rar B'$ where $X'$ is a
desingularization 
 of a flattening of $X$, where any exceptional divisor for $X'/B'$ is
exceptional for $X'/X$. Since $\omega_{B'/B}$ is effective, one immediately
obtains:

\begin{th}[Koll\'ar-Viehweg]\label{kv}  Suppose that $\pi:X\rar B$ is as
above. Then $\omega_\pi$ is big. \qed
\end{th}

\section{PLURI-NODAL REDUCTION}\subtitle{We Prove a Pluri-Nodal Reduction
Lemma and Show That Pluri-Nodal Families Are Mild}

\subsection{Statement} Let $X_0\rar B_0$ be a fibration.  We need to dominate
it by a pluri-nodal fibration, so that it becomes a quotient by the action of a
finite group.  

To this end, we prove the following theorem, which is a variant of de Jong's
results in \cite{dj}, sections 6-7. The proof is based on that of de Jong.

\begin{lem}[{Galois pluri-nodal reduction lemma}]\label{pluri-nodal}
There exists a diagram
$$ \begin{array}{ccc} Y & \rar & X_0 \\
                    \dar &  & \dar \\
                     B_1 & \rar & B_0\end{array} ,$$
and a finite group $G \subset \Aut_{X_0\ttimes_{B_0} B_1}Y$ such that $B_1\rar
B_0$ 
is an 
alteration, $Y/G \rar X_0\ttimes_{B_0} B_1$  is birational  and $Y\rar B_1$ is
a pluri-nodal fibration. 
\end{lem}

{\bf Proof.}
We proceed by induction. The setup is as follows: suppose we have $X\rar Z
\rar B$ a pair of  
fibrations, where $X\rar Z$ is pluri-nodal, and a finite
group $G_0\subset\Aut_B(X\rar Z)$. We also assume that we have a birational
morphism 
$X/G_0\rar X_0\ttimes_{B_0}B$. We will produce a diagram 
$$\begin{array}{ccccccc} X' & \rar &   Z'' & \rar & Z' & \rar    &  B'  \\
                       \dar &      &  \dar &      &    &         & \dar \\
                         X  & \rar &   Z   & & \longrightarrow & &  B 
  \end{array}
$$
with the following properties:
\begin{enumerate} 
  \item the vertical arrows are alterations,
  \item the horizontal arrows are fibrations,
  \item the morphism  $Z''\rar Z'$ is a nodal fibration,
  \item $X' = X\times_Z Z''$, and therefore $X'\rar Z'$ is pluri-nodal,
  \item there is a finite group  
     $G' = G_0\times G'' \subset \Aut_{B'}(X'\rar Z''\rar Z')$, and
  \item the morphism $X'/G'' \rar  X\ttimes_B B'$ is birational, and therefore
       $X'/G'\rar X_0\ttimes_{B_0}B'$ is birational.
\end{enumerate}
 The basis of the induction is $X_0\rar X_0\rar B_0$ with $G_0$ trivial. The
induction ends with $Z'\rar B'$ being birational, in which case we set $Y :=
X', \quad B_1 := Z', \quad G := G'$ and the lemma will be proved.

 Let 
$G_Z\subset\Aut_BZ$ be the image of $G_0$, and denote $W=Z/G_Z$.

\begin{lem} There exists a dominant rational map
$Z/G_Z\das P\rar B$, where $P\rar B$ is a projective bundle, such
that $\dim(P) = \dim(Z)-1$, and such that the generic fiber of $Z$ over $P$ is
geometrically irreducible. 
\end{lem} 
{\bf Proof.} This is obvious in case $\reld(Z/B)=1$, so assume
$\reld(Z/B)>
1$. Denote this relative dimension by $r$. Since we are looking for a rational
map, we may replace $B$ by its 
generic point $\eta$, and replace $Z$ by $Z_\eta$. Let $W=Z/G$, choose an
embedding $W\subset \bfp^N$, and let $f:Z\rar\bfp^N$ be the
induced morphism. According to \cite{jou}, 6.3(4), for general hyperplane
$H\subset\bfp^N$ we have $f^{-1}H$ geometrically irreducible. Continuing by
induction, there is a linear series $(\bfp^{r-1})^*$ of dimension $r-1$ of
hyperplanes in $\bfp^N$ such that the general fiber  of
$Z\das \bfp^{r-1}$ is a geometrically irreducible curve. \qed


The normalization of the closure of the graph of the
rational map $Z\das P$ 
gives a $G_Z$-equivariant resolution of
indeterminacies 
$$ \begin{array}{ccc} Z_1 & \rar & P \\
                    \dar &  &  \\
                     Z &  & \end{array}. $$
 Let $X_1 = X\times_Z Z_1$. Then $X_1\rar Z_1$ is pluri-nodal, and the action
of $G_0$ on $X$ lifts to $X_1$ (if $x_1 = (x,z_1)\in X_1$ and $g\in G_0$ then
$(g(x),g(z_1))\in X_1$ as well). 

We will now perform a canonical nodal reduction for $Z_1\rar P$ using the
Kontsevich space of stable maps.
 The generic 
fiber of $Z_1 \rar P$ is a normal curve, and therefore smooth. Choose a
projective embedding  $Z_1\subset \bfp^N$. Let $d$ be the degree of the generic
fiber of 
$Z_1\rar P$ and let $g$ be its
genus. By 
\cite{bm}, theorem 3.14, there exists a proper Deligne-Mumford
stack  $\OM_{g,0}(Z_1,d)$ parametrizing stable maps $C\rar Z_1$ of
curves of genus $g$ and degree $d$. By \cite{pan} this stack admits a
projective coarse
moduli space. In particular, this implies that there is a finite cover $\rho:M
\rar 
{\OM_{g,0}}(Z_1,d)$ where $M$ is a projective scheme
admitting a stable map $(C\rar M, f:C\rar Z_1)$ whose moduli map is $\rho$. 

Let $\eta\in P$ be the generic point. The pair $((Z_1)_{\eta}\rar \eta,
(Z_1)_{\eta}\hookrightarrow Z_1)$ is a stable map of genus $g$ and degree $d$,
therefore we have a rational map $P\das 
{\OM_{g,0}}(Z_1,d)$. 
We can choose a normal resolution of indeterminacies  
$$ \begin{array}{ccc} P_2 & {\rar} & M\\
                    \dar &  &  \\
                     P &   &  \end{array} $$
such that there is a finite group $G_1\subset\Aut_{P}P_2$ with $P_2/G_1\rar P$
birational.  Let $Z_2 = C\times_MP_2$. We have an induced stable map $(Z_2\rar
P_2,f_2:Z_2\rar Z_1)$, in particular $Z_2\rar P_2$ is nodal. Over the generic
point of $P_2$ this coincides with $Z_1\times_{P} P_2$. 

Since stable reduction
over a normal base is unique when it exists (see \cite{d-o}, 2.3), the action
of $G_1$ lifts to 
$Z_2$, and  it lifts  
to $X_2 = X_1\times_{Z_1}Z_2$ as well by pulling back as before. Let $P_2\rar
B_2\rar B$ be the Stein 
factorization. Since the Stein factorization is unique we have canonically an
action of $G_1$ on $B_2$. Let $G_2\subset G_1$ be the subgroup acting trivially
on 
$B_2$. Then 
$G=G_0\times G_2\subset\Aut_{B_2}(X_2\rar P_2)$. We have $X_2\rar P_2$
pluri-nodal, 
and 
$X_2/G_2\rar X\ttimes_B B_2$ birational. If we denote $X' := X_2, \quad Z''
:= 
Z_2, \quad Z' :=
P_2, \quad B' := B_2 $ and $G'' := G_2$ we have obtained the goal of the
induction step. \qed

\subsection{Mild Singularities}\label{mild}

We want to show that pluri-nodal fibrations are mild. This seems to be well
known (see \cite{hassett}, \S 4), but in our case we can give a  
proof which is sufficiently short to include here. The following lemma is
well known (see \cite{viehweg},  lemma 3.6): 
\begin{lem}
Let $Y\rar B$ be a nodal fibration such that $B$ is smooth and the discriminant
locus is a divisor of normal crossings. Then $Y$ is r-G.
\end{lem}
{
(The proof is by taking formal coordinates near a singular point of the form
$xy = t_1^{k_1}\cdots t_r^{k^r}$, and either resolving singularities
explicitly or noting that this is a toroidal singularity.)}

\begin{prp} 
Let $Y\rar B$ be a nodal fibration such that $B$ is r-G. Then $Y$ is r-G.
\end{prp}

{\bf Proof.} Let $r:B_1 \rar B$ be a resolution of singularities, $Y_1\rar B_1$
the pullback, and assume that the discriminant locus of $Y_1\rar B_1$ is a
divisor of normal crossings. Let $f:Y_1\rar Y$ be the induced map. Then
$r_*\omega_{B_1}  = \omega_B$ and $f^*\omega_{Y/B} = \omega_{Y_1/B_1}$, and by
the projection formula we obtain that $f_*\omega_{Y_1} = \omega_Y$.\qed

By induction we obtain:
\begin{cor}
If $\pi:Y\rar B$ is a pluri-nodal fibration where $B$ is r-G, then $Y$ is r-G.
In particular the $n$-th fibered power $Y^n_B$ is r-G.
\end{cor}

Thus pluri-nodal fibrations are mild.

\section{PROOF OF THE THEOREM}\subtitle{Our Main Theorem Arrives at an
Enchanted Place, and We Leave It There.} 
Let $X_0\rar B_0$ be a smooth projective family 
of varieties 
of general type of maximal variation. Choose a  model
${X}\rar {B}$ where both $X$ and  $B$ are projective nonsingular. By 
\ref{pluri-nodal} we may 
assume, after an 
alteration $B_1\rar {B}$, that we have a birational morphism
$g_0:Y/G=X_1 \rar X\ttimes_{B}B_1$ where
$\pi_Y:Y\rar B_1$ is a pluri-nodal 
fibration and $G\subset\Aut_{B_1}Y_1$ a finite group. Choose a resolution of
singularities $r:X_2\rar X_1$ and denote by $\pi_2:X_2\rar B_1$ the projection.
We have a diagram:
\begin{equation}\label{diagram}
\begin{array}{lclcl}
 	&		 & Y    & 		     &	 \\
 	&		 & \dar q & 		     &	 \\
  X_2 & \stackrel{r}{\lrar}& X_1 & \stackrel{g_0}{\lrar} & X \\
     &  \searrow^{\pi_2}  &  \dar &                    &   \dar \\
     &                  & B_1   &   \rar             & B 
\end{array}
\end{equation}
 According to
\ref{product} (where we set $V=Y$ and $W=X_1$) there is an ideal $\I_F\subset
\co_{B_1}$  
such that $(\I_F^{e(G)\cdot n})\cdot\J_{m,n} \subset \J_n(G^m,Y^m_{B_1})$.
For arbitrary integer $m>0$ let $\clx_m\rar X^m_B$ be a resolution of
singularities of the main component, and let $\W_m\rar (X_1)^m_{B_1}$ be a
resolution of 
singularities, dominating $\clx_m$. 
According to \ref{base-change} (applied to $B_1\rar B$) there is an ideal
$\I\subset \co_{B_1}$, such that for any $m$, if $\omega_{\W_m/B_1}\otimes \I$
is big then $\clx_m$ (and therefore $(X_0)^m_{B_0}$) is of general type.

By the Koll\'ar - Viehweg theorem \ref{kv}, $\omega_{\pi_2}$ is big. Therefore
$q^*r_*\omega_{\pi_2}$ is big. We
have by definition that $\omega_{\pi_Y}\otimes \J_1(G,Y)$ is big.
 Therefore,
for sufficiently large $n$ we have that $\omega^{n}_{\pi_Y}\otimes
{\J_n}\I\I_F^{e(G)}$ is big. Pulling back along all the projections
$p_i:Y^m_{B_1}\rar 
Y$ we  have that $\omega^n_{Y^m_{B_1}/B_1}\otimes {\J_{m,n}}\I^m\I_F^{m\cdot
e(G)}$  
is big. In particular, if $m>n$, we have that
$\omega^n_{Y^m_{B_1}/{B_1}}\otimes {\J_{m,n}}\I^n\I_F^{n\cdot e(G)}$ is big. By
\ref{product} we have that $\omega^n_{Y^m_{B_1}/{B_1}}\otimes
\J_n(G^m,Y^m_{B_1})\I^n$ is big.
Taking invariants and using \ref{invariants}, $\omega_{\W_m/B_1}\otimes \I$ is
big,  and by \ref{base-change} we have that $(X_0)^m_{B_0}$
is of general type for large $m$. 

\subsection{An alternative approach.}\label{alternative} The following argument
gives a 
variation on the proof which is more in line with \cite{kollar} and
\cite{viehweg}. Having chosen the diagram (\ref{diagram}), we can alter it as
follows: using semistable reduction in codimension 1 (see \cite{te} II, and
\cite{kawamata}, theorem 17), we
can find a nonsingular alteration  $B_1'\rar B_1$, a 
variety $X_2'\rar B_1'$, and a birational morphism $X_2'\rar
X_2\ttimes_{B_1}B_1'$  satisfying the following conditions
\begin{enumerate} 
\item The discriminant locus $\Delta$ of $X_2'\rar B_1'$ is a divisor of normal
crossings. Set $F=Sing(\Delta)$ and $U=B_1'\setminus F
\stackrel{i}{\hookrightarrow} B_1'$.
\item  The restriction $X_2'|_U\rar U$ is semistable, in particular it is mild
(see \cite{viehweg}, lemma 3.6).
\end{enumerate}
 Let $X_1'= X_1\times_{B_1}B_1'$
and $Y'= Y\times_{B_1}B_1'$. We can replace $B_1,X_1,X_2, Y$ by
$B_1',X_1',X_2', Y'$ and assume that conditions (1) and (2) are satisfied.

 Let $\pi_{X_{(m)}}:X_{(m)}\rar B_1$ be the main component of
$(X_1)^m_B$. Choose a resolution of singularities $W_m\rar X_{(m)} $, and
let 
$\pi_{W_m}:W_m \rar B_1$ be the associated projection. Denote  $\F_{m,n}=
{\pi_2}_*\omega_{\pi_{W_m}}$ and    $\G_{m,n}=(\F_{m,n})^{**}$. Since the
restriction of $W_m$ to $U$ is mild, we have that  $\G_{m,n} =
i_*i^*\F_{m,n}$. Applying \ref{smoothproduct}, we obtain:
\begin{enumerate}
\item We have natural morphisms 
  $$ \G_{m_1,n}\otimes  \G_{m_2,n} \rar \G_{m_1+m_2,n} $$
(by pulling back sections to $W_{m_1+m_2}$ over $U$, multiplying and
extending). 
\item We have natural morphisms $$ \G_{m,n_1}\otimes  \G_{m,n_2} \rar
\G_{m,n_1+n_2} $$ 
(by multiplying sections).
\item We have $$\G_{m,n}\otimes \I_F^{n(e(G)-1)} \subset \F_{m,n}\subset
\G_{m,n}$$ (by \ref{product}. Notice that the remark after \ref{product} shows
that  $\I_F^n$  suffices). 
\end{enumerate}

By Kollar's theorem $G_{1,n}$ is \bbig\ for sufficiently large $n$.  We can
choose an ideal $\I$ as in \ref{base-change}. By (2) above, for sufficiently
large $n$ we have that  $G_{1,n}\otimes \I\I_F^{e(G)}$ is big, and using (1)
above we have that for sufficiently large $m$, we have that  $G_{m,n}\otimes
\I^m\I_F^{n\cdot e(G)}$ is big, therefore by (3)  $\F_{m,n}\otimes\I^m$ is big,
which is what we need.

\if\FUN y
\vspace*{3mm}
\hspace*{2.5in}\parbox{4in}{\em ...``Come on!'' \\ ``Where?'' said Pooh. \\
``Anywhere,'' Said Christopher Robin. \\ \flushright{\em A.A. Milne,\\ The
House at Pooh Corner}}
\fi
\if\figure y
\vspace*{2.in}
\hspace*{-1in}
\flushleft{
\begin{picture}(5,5)
\put(-1,0){\psfig{file=pooh1.ps}}
\end{picture}}
\fi


\begin{thebibliography}{HHHHHHH}

 \bibitem[$\aleph$]{a} D. Abramovich, {\em Uniformit\'e des points rationnels
des courbes alg\'ebriques sur les extensions quadratiques et cubiques},
 C.R. Acad. Sc. Paris, t. 321, S\'er. I, p. 755-758, 1995.

 \bibitem[$\aleph$1]{a-int}  D. Abramovich, {\em Uniformity of stably integral
points on elliptic curves,} preprint.\\
{\tt http://math.bu.edu/INDIVIDUAL/abrmovic/integral.ps}

\bibitem[$\aleph$-V]{av} D. Abramovich and J. F. Voloch, {\em Lang's
conjecture, fibered powers and uniformity}, preprint.\\
{\tt http://math.bu.edu/INDIVIDUAL/abrmovic/conjh.dvi}

\bibitem[B-M]{bm} K. Behrend and Yu. Manin, {\em Stacks of Stable Maps and
Gromov-Witten Invariants}, preprint.\\
{\tt http://xxx.lanl.gov/e-print/alg-geom/9506023}

\bibitem[CHM]{chm}  L. Caporaso, J. Harris, B. Mazur: {\em Uniformity of 
   rational points,}  J. Amer. Math. Soc., to appear. \\
{\tt ftp://ftp.math.harvard.edu/pub/uniformityofrationalpoints.tex}

\bibitem[D-M]{d-m} M. Martin - Deschamps and R. Lewin - Menegaux, {\em
Applications rationnelles s\'eparables dominantes sur une vari\'et\'e de type
g\'en\'eral, } Bull. S. M. France 106 (1978) no. 3, p. 279-287.

\bibitem[E]{elkik} R. Elkik, {\em   Singularit\'es rationnelles
et d\'eformations},  Inv.  Math. 47, 1978, p.
139-147.


 \bibitem[Has]{hassett} B. Hassett: {\em Correlation for surfaces of general
type}, preprint, 1995.\\
{\tt   http://xxx.lanl.gov/e-print/alg-geom/9507015}

\bibitem[dJ]{dj} A. J. de Jong, {\em Smoothness, semistability, and
alterations}, preprint.\\ 
{\tt ftp://ftp.math.harvard.edu/pub/AJdeJong/alterations.dvi}

\bibitem[dJ-O]{d-o} A. J. de Jong and F. Oort, {\em On extending families of
curves}, preprint.\\  
{\tt ftp://ftp.math.harvard.edu/pub/AJdeJong/extend.dvi}

\bibitem[J]{jou} J.-P. Jouanolou, {\em Th\'eor\`emes de Bertini et
Applications}, Progress in Math. 42, Birkh\"auser: Boston, Basel, Stuttgart,
1983. 

\bibitem[Ka]{kawamata} Y. Kawamata, {\em Characterization of abelian
varieties,} Comp. Math. 43, 1981 p 253-276.

\bibitem[KKMS]{te} G. Kempf, F. Knudsen, D. Mumford and B. Saint-Donat,
{\em Toroidal Embeddings I}, Springer, LNM 339, 1973.

\bibitem[Ko]{kollar} J. Koll\'ar, {\em Subadditivity of the Kodaira dimension:
   fibers of general type}.  Adv. Stud. in Pure Math. 10 (1987)
	p. 361-398.

\bibitem[L]{langbul}  S. Lang, {\it Hyperbolic diophantine analysis.}
   Bull. A.M.S. 14 (1986) p. 159-205.

\bibitem[Mor]{moriwaki} A. Moriwaki, {\em Remarks on S. Lang's conjecture over
function fields,} preprint.\\
{\tt   http://xxx.lanl.gov/e-print/alg-geom/9412021}

\bibitem[Pac]{p} P. Pacelli,  {\em Uniform boundedness for rational points},
preprint. 

\bibitem[Pan]{pan} R. Pandharipande, {\em Notes On Kontsevich's
Compactification Of The Space Of Maps}, preprint.

\bibitem[V1]{viehweg1} E. Viehweg,  {\em Die Additivit\"at der Kodaira
Dimension f\"ur projektive Fasserr\"aume \"uber Variet\"aten des allgemeinen
Typs,} Jour. reine und angew. Math. 330 (1982), 132-142.

\bibitem[V2]{viehweg} E. Viehweg, {\em Weak positivity and the additivity of
the Kodaira dimension for certain fiber spaces.} In: Algebraic varieties and
analytic varieties, Advanced Studies in Pure Math. 1 (1983), 329 - 353.

\bibitem[V3]{vletter} E. Viehweg, letter to Barry Mazur.

\end{thebibliography}
\end{document}